\documentclass[11pt]{article}
\usepackage{amsfonts}
\usepackage{amsmath}
\usepackage{amssymb}
\usepackage[margin=1in]{geometry}

\newtheorem{theorem}{Theorem}
\newtheorem{lemma}[theorem]{Lemma}
\newtheorem{proposition}[theorem]{Proposition}
\newtheorem{example}[theorem]{Example}
\newtheorem{remark}[theorem]{Remark}
\newenvironment{proof}[1][Proof]{\noindent\textbf{#1.} }{\ \rule{0.5em}{0.5em}}

\begin{document}

\title{Zero extension for the Dirichlet problem on graphs}
\author{Yanglan Ou\thanks{Corresponding authors. ownerounan@163.com (Y. Ou).} \\
School of Pharmacy, Guangzhou Xinhua University\\
Guangzhou 510520, Guangdong, P. R. China \and Xuechao Yao \\
School of Integrated Circuits and New Energy\\
Guangzhou College of Technology and Business\\
Guangzhou 510850, Guangdong, China \and Guang Zhang \\
Zhujiang College, South China Agricultural University\\
Guangzhou 510900, Guangdong, P. R. China}
\date{}
\maketitle

\begin{abstract}
Solutions of the Dirichlet problem on a finite part of a weighted graph are
extended by zero, together with their source, to a larger connected part. A
necessary and sufficient condition is proved for the extended function to
remain a solution of the corresponding Dirichlet problem on the enlarged
domain: at every newly included vertex adjacent to the original part, the
weighted sum of the solution over its neighbors in that part must vanish.
The admissible sources are characterized as the orthogonal complement of the
harmonic extensions from the enlarged domain to the original part, their
space has dimension equal to the number of vertices of the original part
minus the rank of the interface coupling, and an exact energy identity for
the extension error follows from a Schur complement of the graph Laplacian.
The same reduction covers other boundary conditions, including the Neumann
problem, and extends to higher powers of the Laplacian, where the boundary
data must be compatible with the lower-order interface traces. The discrete
criterion is compared with its continuous counterpart, and no result in this
form has been reported on Riemannian manifolds.
\end{abstract}



\section{Introduction}

Throughout, $\Gamma=(V,E,w)$ is a connected locally finite weighted graph,
and $\mathcal{L}$ denotes its combinatorial Laplacian,
\begin{equation*}
\mathcal{L}u(x)=\sum_{y\sim x}w_{x,y}\bigl(u(y)-u(x)\bigr).
\end{equation*}
Fix a nonempty finite connected set $F\subset V$, with exterior boundary $%
\delta(F)=\{x\in V\backslash F:x\sim y\text{ for some }y\in F\}$ and closure
$\overline{F}=F\cup\delta(F)$, and take a source $f\in\mathbb{R}^{F}$. The
Dirichlet problem on $F$,
\begin{equation}  \label{eq:dir}
\begin{cases}
-\mathcal{L}u(x)=f(x), & x\in F, \\
u(x)=0, & x\in\delta(F),%
\end{cases}%
\end{equation}
has a unique solution $u$ by the maximum principle. Now enlarge the domain:
take a nonempty finite connected set $\Omega$ with $F\subseteq\Omega%
\subseteq V$ and $\Omega^{c}\neq\emptyset$, write $S=\Omega\backslash F$,
and extend $u$ and $f$ by zero to $\Omega$,
\begin{equation}  \label{eq:ext}
\widetilde{u}(x)=%
\begin{cases}
u(x), & x\in F, \\
0, & x\in\overline{\Omega}\backslash F,%
\end{cases}
\qquad \widetilde{f}(x)=%
\begin{cases}
f(x), & x\in F, \\
0, & x\in\overline{\Omega}\backslash F.%
\end{cases}%
\end{equation}
The zero extension problem asks whether $\widetilde{u}$ solves the same
Dirichlet problem on the larger domain,
\begin{equation}  \label{eq:extprob}
\begin{cases}
-\mathcal{L}v(x)=\widetilde{f}(x), & x\in\Omega, \\
v(x)=0, & x\in\delta(\Omega).%
\end{cases}%
\end{equation}
Since \eqref{eq:ext} preserves the equation on $F$ by construction, the
question concerns only the newly included vertices $S$: it asks when the
weighted couplings from $S$ to $F$ produce no residual, and, if they do, how
large the resulting error is.

This is the problem studied in the paper. We prove that \eqref{eq:extprob}
holds if and only if the weighted sum of $u$ over the neighbors in $F$
vanishes at every interface vertex, that is, if and only if the interface
block of $\mathcal{L}$ annihilates $u$ on $(\Omega \backslash F)\cap \delta
(F)$; no condition is imposed on the exterior boundary $\delta (\Omega )$,
and a four-vertex example shows that imposing one is incorrect. We then
characterize the admissible sources by orthogonality to the harmonic
extensions from $S$ to $F$, give the dimension of their space, and show that
for independently varying positive interface weights its generic dimension
equals $|F|$ minus the maximum matching size of the interface bipartite
graph. When compatibility fails, the Schur complement of $\mathcal{L}$ on $%
\Omega $ with respect to $F$ gives the exact energy of the extension error,
together with two-sided spectral bounds; a weight-perturbation example shows
that exact compatibility can disappear under an arbitrarily small change of
the weights while the solution error still tends to zero.

For continuum Poisson equations, Cai and Zhou \cite{Cai} studied when zero
extension preserves the equation in classical, strong and weak formulations.
On a graph the question has an elementary answer because the Laplacian has
finite range, but the exact criterion above leaves further questions, namely
the size of the admissible-source space, its dependence on the interface
topology and weights, and the behaviour when compatibility holds only
approximately. These are answered by the results below. The constructions we
use are standard: the map from interior sources to interface residuals is
the accompanying matrix in Kron reduction, and Schur complements with their
network interpretations are treated in \cite{Dorfler}; the generic rank of
an interface block is governed by bipartite matching \cite{Davis}. We use
these constructions explicitly rather than presenting elementary block
identities as a new theory. The interface reduction itself is not specific
to Dirichlet data: other boundary value problems on networks can be recast
in the same form through equilibrium measures on the boundary \cite{Bend}.
All results concern a fixed graph or a finite family of weights; no
continuum-limit or empirical imaging claim is made. On the continuous side,
the zero extension criterion for the Poisson equation is known \cite{Cai},
and for the biharmonic equation it has been settled by Xu and Zhou \cite{Xuz}
and Xu and Yu \cite{xuy}; to the best of our knowledge, however, no result
in the form developed here---an interface criterion with an exact error
identity, together with its higher-order compatibility constraints---has
been reported on Riemannian manifolds, where the corresponding statements
are phrased through normal derivatives, harmonic orthogonality, or unique
continuation rather than through a finite interface operator.

The criterion is stated for Dirichlet data, but the same interface reduction
covers other boundary conditions: boundary value problems on weighted
networks, including Neumann and mixed problems, can be recast in this form
through equilibrium measures on the boundary \cite{Bend}, which in the
present notation amounts to prescribing the interface trace of the extension
rather than setting it to zero. Thus the Neumann problem is included as a
particular case, with the interface flux $\mathcal{L}_{SF}u$ itself, rather
than its vanishing, as the prescribed quantity.

For $n\geq 1$ define $\mathcal{L}^{n}$ recursively by $\mathcal{L}^{1}=%
\mathcal{L}$ and $\mathcal{L}^{n}=\mathcal{L}\circ \mathcal{L}^{n-1}$, and
consider $(-\mathcal{L})^{n}u=f$ on $F$. For $n\geq 2$ the restricted
operator and the restriction of the power no longer coincide: block
multiplication on $\Omega =F\cup S$ gives $E^{T}\mathcal{L}_{\Omega }^{2}E=%
\mathcal{L}_{F}^{2}+\mathcal{L}_{FS}\mathcal{L}_{SF}$, and more generally $%
E^{T}\mathcal{L}_{\Omega }^{n}E=\mathcal{L}_{F}^{n}+$ terms supported on the
interface, so the higher-order analogue of \eqref{eq:cond} is not obtained
by replacing $\mathcal{L}_{F}$ with $\mathcal{L}_{F}^{n}$ in Theorem~\ref%
{thm:crit}. The boundary data for $(-\mathcal{L})^{n}$ are not free: they
must be compatible with the lower-order interface traces already fixed on $F$%
. For $n=2$, the continuous biharmonic Dirichlet problem requires $u=0$ and $%
\partial _{\nu }u=0$ on $\delta (F)$, and on a graph the second condition is
replaced by a prescription of the interface trace of $\mathcal{L}u$; since $%
\mathcal{L}u$ on $\delta (F)$ is determined by $u$ and its neighbors,
specifying the two traces independently overdetermines the problem unless
the pair lies in the range of the interface operator. The interface block
persists at every order, but the analogue of $\mathcal{C}$ for $(-\mathcal{L}%
)^{n}$ is the kernel of a block interface operator that encodes these
compatibility constraints among the iterated traces, rather than of $%
\mathcal{L}_{SF}$ alone. The continuous biharmonic criterion is known \cite%
{Xuz,xuy}; for $n\geq 3$ on graphs, and for the corresponding continuous
higher-order interface criteria, we have not found such results in the
literature.

The criterion has a direct network reading and a range of candidate
applications. In a grounded resistor network, $u$ is a voltage, $f$ an
injected current, and the interface residual the net current required at the
newly included zero-voltage vertices. The same structure appears when a
quantity known on a subregion is continued by zero and one asks whether the
continuation still satisfies the governing equation: graph-based image
inpainting and point-cloud denoising \cite{Bend}, missing-data infilling in
zero echo-time and undersampled MRI and in photon-counting CT \cite%
{Huo,Lian,Chen}, and Laplacian-regularized matrix factorization for
drug--target interaction prediction \cite{mon,Wangm}. In each case the
present paper supplies only the compatibility criterion and the error
identity, not a validated algorithm; establishing the latter would require
its own objective, graph construction and numerical study.

\section{Main results}

\label{sec:main}

For $K\subseteq V$ let $\mathcal{L}_K$ denote the Dirichlet restriction of $-%
\mathcal{L}$ to $K$, that is, the operator acting on functions vanishing on $%
\delta(K)$ and defined by $-\mathcal{L}$ at the vertices of $K$; its matrix
has diagonal entries $\sum_{y\in V}w_{x,y}$ and off-diagonal entries $%
-w_{x,y}$ for $x\neq y\in K$.

\begin{lemma}
\label{lem:pd} For every nonempty finite proper subset $K\subset V$, the
operator $\mathcal{L}_K$ is symmetric and positive definite. Hence %
\eqref{eq:dir} has the unique solution $u=\mathcal{L}_F^{-1}f$ on $F$.
\end{lemma}

\begin{proof}
For $z\in\mathbb{R}^{K}$, symmetry gives $z^{T}\mathcal{L}_K
z=\tfrac12\sum_{x,y}w_{x,y}(E_{K}z(x)-E_{K}z(y))^{2}\ge0$, where $E_{K}z$ is
the extension by zero. If the sum vanishes, $E_{K}z$ is constant along every
edge, hence on the connected graph, and it is zero on the nonempty set $%
V\backslash K$.
\end{proof}

\bigskip

At the level of $\mathcal{L}$, the operator on $\Omega $ decomposes as
\begin{equation}
\mathcal{L}_{\Omega }=%
\begin{pmatrix}
\mathcal{L}_{F} & -\mathcal{L}_{FS} \\
-\mathcal{L}_{SF} & \mathcal{L}_{S}%
\end{pmatrix}%
,  \label{eq:block}
\end{equation}%
where $\mathcal{L}_{FS}$ and $\mathcal{L}_{SF}$ are the interface blocks,
with entries $w_{x,y}$ for $x\in F$, $y\in S$, and $\mathcal{L}_{S}$ is the
Dirichlet restriction to $S$ with the couplings to $F$ included in its
diagonal. The block $\mathcal{L}_{SF}$ vanishes on rows indexed by vertices
of $S$ that are not adjacent to $F$, so the only nontrivial equations below
are interface equations.

\begin{theorem}[Interface criterion]
\label{thm:crit} The extension $\widetilde{u}$ defined by \eqref{eq:ext}
solves \eqref{eq:extprob} if and only if
\begin{equation}  \label{eq:cond}
\sum_{y\in F,\,y\sim x}w_{x,y}u(y)=0 \qquad\text{for every }%
x\in(\Omega\backslash F)\cap\delta(F),
\end{equation}
that is, if and only if $\mathcal{L}_{SF}u$ vanishes on the interface. For $%
x\in(\Omega\backslash F)\backslash\delta(F)$ the sum in \eqref{eq:cond} is
empty, so \eqref{eq:cond} is equivalent to the same identity for every $%
x\in\Omega\backslash F$. In that case $\mathcal{L}_\Omega E u-E\widetilde{f}%
=(0,-\mathcal{L}_{SF}\mathcal{L}_F^{-1}f)^{T}$, where $Eu=(u,0)$ and $E%
\widetilde{f}=(f,0)$.
\end{theorem}

\begin{proof}
At $x\in F$, every neighbor outside $F$ has zero value in both $u$ and $%
\widetilde{u}$, so $-\mathcal{L}\widetilde{u}(x)=f(x)$. At $x\in S$, $%
\widetilde{u}(x)=0$, hence $-\mathcal{L}\widetilde{u}(x)=-\sum_{y\in
F,\,y\sim x}w_{x,y}u(y)=-(\mathcal{L}_{SF}u)_{x}$. The right-hand side $%
\widetilde{f}$ is zero on $S$, and the values on $\delta(\Omega)$ already
vanish.
\end{proof}

\bigskip

The condition \eqref{eq:cond} is imposed only at the vertices of $(\Omega
\backslash F)\cap \delta (F)$. Vertices of $\delta (F)\backslash \Omega $
are adjacent to $F$ but lie outside the enlarged domain; \eqref{eq:extprob}
imposes no equation there, and \eqref{eq:cond} is silent about them.
Vertices of $\delta (F)\cap \delta (\Omega )$ are adjacent to $F$ and lie on
the exterior boundary of $\Omega $; there \eqref{eq:extprob} prescribes only
the Dirichlet value $v=0$, so a nonzero interface residual is allowed, as
Example~\ref{ex:counter} illustrates. If $\delta (F)\subset \Omega $, every
vertex adjacent to $F$ belongs to $\Omega \backslash F$, and \eqref{eq:cond}
covers the whole interface; this is the case in which the condition can be
read as a complete interface condition. No equation is imposed on $\delta
(\Omega )$ in any case, and replacing $(\Omega \backslash F)\cap \delta (F)$
in \eqref{eq:cond} by $(S\cup \delta (\Omega ))\cap \delta (F)$ gives a
sufficient condition but, in general, not a necessary one.

\begin{example}[Exterior boundary counterexample]
\label{ex:counter} Take vertices $a,b,c,d$ with unit-weight edges $%
ab,ac,bc,ad$, and set $F=\{a,b\}$, $\Omega=\{a,b,c\}$. For $u=(1,-1)^{T}$, $%
f=\mathcal{L}_F u=(4,-3)^{T}$ and $\mathcal{L}%
_\Omega(1,-1,0)^{T}=(4,-3,0)^{T}$, so the extension is the required
solution. Nevertheless $d\in\delta(\Omega)$ and $\sum_{y\in F}w_{d,y}u(y)=1$;
imposing an equation at $d$ would incorrectly exclude this solution.
\end{example}

\bigskip

Define the compatible-source space by
\begin{equation}
\mathcal{C}=\{f\in \mathbb{R}^{F}:E\mathcal{L}_{F}^{-1}f=\mathcal{L}_{\Omega
}^{-1}E\widetilde{f}\}.  \label{eq:C}
\end{equation}%
For $g\in \mathbb{R}^{S}$ let $H_{g}=\mathcal{L}_{F}^{-1}\mathcal{L}_{FS}g$;
this is the restriction to $F$ of the function harmonic on $F$ with
prescribed values $g$ on $S$ and zero on $V\backslash \Omega $.

\begin{theorem}[Dual description]
\label{thm:dual} $\mathcal{C}=\ker(\mathcal{L}_{SF}\mathcal{L}_F^{-1})=%
\mathcal{L}_F(\ker\mathcal{L}_{SF})=(\operatorname{ran}(\mathcal{L}_F^{-1}\mathcal{L}%
_{FS}))^{\perp}$, and $\dim\mathcal{C}=|F|-\operatorname{rank}\mathcal{L}_{SF}$.
Equivalently, $f$ is compatible if and only if
\begin{equation}  \label{eq:orth}
\sum_{y\in F}f(y)H_{g}(y)=0\qquad\text{for every }g\in\mathbb{R}^{S}.
\end{equation}
\end{theorem}

\begin{proof}
Theorem~\ref{thm:crit} gives $\mathcal{C}=\ker(\mathcal{L}_{SF}\mathcal{L}%
_F^{-1})=\mathcal{L}_F(\ker\mathcal{L}_{SF})$. Since $\mathcal{L}_F^{-1}$ is
symmetric, $g^{T}\mathcal{L}_{SF}\mathcal{L}_F^{-1}f=f^{T}\mathcal{L}_F^{-1}%
\mathcal{L}_{FS}g=f^{T}H_{g}$. Multiplication by $\mathcal{L}_F^{-1}$
preserves rank, so rank-nullity gives the dimension formula.
\end{proof}

\bigskip

The criterion is local in $u$ but generally nonlocal in $f$ because $%
\mathcal{L}_{F}^{-1}$ need not be sparse. If $\delta (F)\subset \Omega $, $g$
can prescribe arbitrary values on the whole exterior boundary of $F$, and %
\eqref{eq:orth} tests all harmonic extensions to $F$. If $F\subseteq \Omega
_{1}\subset \Omega _{2}$, the compatible-source space for $\Omega _{2}$ is
contained in that for $\Omega _{1}$, and once $\Omega $ contains $\delta (F)$
further enlargement does not change it.

Consider the bipartite interface graph $\mathcal{I}$ with vertex classes $F$
and $S$ and an edge $y\sim x$ whenever $w_{x,y}>0$, and let $\nu(\mathcal{I}%
) $ be its maximum matching size.

\begin{proposition}[Generic dimension]
\label{prop:gen} Fix the graph topology and all non-interface weights, and
allow the positive weights on distinct edges of $\mathcal{I}$ to vary
independently. Then $\dim\mathcal{C}\ge|F|-\nu(\mathcal{I})$ for every
choice of weights, with equality outside a proper algebraic subset of the
positive weight space.
\end{proposition}

\begin{proof}[Proof (sketch)]
Every nonzero term in a $k\times k$ minor of $\mathcal{L}_{SF}$ selects $k$
entries in different rows and columns, hence a matching, so $\operatorname{rank}%
\mathcal{L}_{SF}\le\nu(\mathcal{I})$. Conversely, a maximum matching yields
a square submatrix whose determinant is a nonzero polynomial in the
independent interface weights \cite{Davis}; its zero set is proper, and
outside it $\operatorname{rank}\mathcal{L}_{SF}\ge\nu$. Apply Theorem~\ref{thm:dual}.
\end{proof}

\bigskip

Independence of the weights matters: equal-weight constraints or geometric
parameterizations can force the weights into an exceptional algebraic set. A
nonzero source $f\geq 0$ or $f\leq 0$ is never compatible: the maximum
principle gives $u>0$ on $F$, and since $\Omega $ is connected and properly
contains $F$, some interface entry of $\mathcal{L}_{SF}u$ is strictly
positive. Thus an exactly localized response requires a source with both
signs.

Let $\nu=\mathcal{L}_\Omega^{-1}E\widetilde{f}$ be the true enlarged-domain
solution, $e=\nu-Eu$, and $r=\mathcal{L}_{SF}u$. The Schur complement of $%
\mathcal{L}_\Omega$ with respect to $F$ is
\begin{equation}  \label{eq:Q}
\mathcal{Q}=\mathcal{L}_S-\mathcal{L}_{SF}\mathcal{L}_F^{-1}\mathcal{L}_{FS}.
\end{equation}

\begin{theorem}[Exact error and energy]
\label{thm:err} The operator $\mathcal{Q}$ is symmetric positive definite,
and
\begin{equation}  \label{eq:err}
e_S=\mathcal{Q}^{-1}r,\qquad e_F=\mathcal{L}_F^{-1}\mathcal{L}_{FS}\mathcal{Q%
}^{-1}r,\qquad \|e\|_{\mathcal{L}_\Omega}^{2}=r^{T}\mathcal{Q}^{-1}r,
\end{equation}
with $\|e\|_{\mathcal{L}_\Omega}^{2}=e^{T}\mathcal{L}_\Omega e$.
Consequently
\begin{equation}  \label{eq:bound}
\frac{\|r\|_{2}^{2}}{\lambda_{\max}(\mathcal{Q})}\le\|e\|_{\mathcal{L}%
_\Omega}^{2}\le \frac{\|r\|_{2}^{2}}{\lambda_{\min}(\mathcal{Q})}.
\end{equation}
\end{theorem}

\begin{proof}
For $z\ne0$, the positive quadratic form of $\mathcal{L}_\Omega$ evaluated
at $(\mathcal{L}_F^{-1}\mathcal{L}_{FS}z,z)$ gives $z^{T}\mathcal{Q}z>0$.
Since $\mathcal{L}_\Omega e=(0,r)^{T}$, the first block equation gives $e_F=%
\mathcal{L}_F^{-1}\mathcal{L}_{FS}e_S$ and the second gives $\mathcal{Q}%
e_S=r $; moreover $e^{T}\mathcal{L}_\Omega e=e_S^{T}r=r^{T}\mathcal{Q}^{-1}r$.
\end{proof}

\bigskip

The identity expresses a standard block elimination in terms of the
extension error; the map $f\mapsto \mathcal{L}_{SF}\mathcal{L}_{F}^{-1}f$ is
the source-to-interface map familiar from Kron reduction \cite{Dorfler}. A
small Euclidean residual alone is not a uniform guarantee of small error
over graph families for which $\lambda _{\min }(\mathcal{Q})$ tends to zero.
If adjustment of the source is part of the modeling problem, the nearest
compatible source in the Euclidean norm is the orthogonal projection of $f$
onto $\ker (\mathcal{L}_{SF}\mathcal{L}_{F}^{-1})$; when the source is
fixed, Theorem~\ref{thm:err} instead provides the actual enlarged-domain
solution.

\subsection*{Two examples}

\emph{Weight perturbation.} Let $V=\{a,b,c,d,g\}$, $F=\{a,b\}$, $%
\Omega=\{a,b,c,d\}$, with unit-weight edges $ab,ac,cd,gd$ and edge $bd$ of
weight $1+\epsilon$, $\epsilon\ge0$. At $\epsilon=0$, $\mathcal{L}_{SF}$ has
rank one and $u=(1,-1)^{T}$ gives the compatible source $f_0=(4,-4)^{T}$.
For every $\epsilon>0$, $\det\mathcal{L}_{SF}=\epsilon$, so $\mathcal{C}%
_\epsilon=\{0\}$. Holding the candidate interior solution fixed, choose
$f_\epsilon=\mathcal{L}_F u=(4,-4-\epsilon)^{T}$; then $r_\epsilon=\mathcal{L}%
_{SF}u=(0,-\epsilon)^{T}$ and $\|r_\epsilon\|_2=\epsilon$. Direct evaluation
of the two-dimensional reduced system gives
\begin{equation*}
\|e_\epsilon\|_{\mathcal{L}_{\Omega_\epsilon}}^{2}=\frac{8\epsilon^{2}(%
\epsilon+2)}{13\epsilon+24}=O(\epsilon^{2}),
\end{equation*}
while for every $\epsilon>0$ one has $\operatorname{dist}(f_\epsilon,\mathcal{C}%
_\epsilon)=\|f_\epsilon\|_2\to4\sqrt{2}$. Thus the extension can be accurate
even when the source is far from the exact compatibility space; the
constants in \eqref{eq:bound} degenerate because the smallest nonzero
singular value of $\mathcal{L}_{SF}\mathcal{L}_F^{-1}$ tends to zero. The
interface graph has a matching of size two, so the one-dimensional
compatibility space at equal weights is exceptional, as Proposition~\ref%
{prop:gen} predicts.

\emph{Sign-changing inhomogeneous localization.} On the integer path with
unit weights, take $F=\{1,2,3\}$ and $u=(0,1,0)^{T}$. For any finite
connected $\Omega$ containing $\{0,1,2,3,4\}$, the interface block
annihilates $u$, the source is $f=(-1,2,-1)^{T}$, and zero extension solves
the enlarged problem. This is an inhomogeneous solution with a deliberately
sign-changing source; it is not a nonzero solution of the homogeneous
equation and does not contradict a homogeneous unique continuation theorem.

\bigskip

Finally, it is worth recording how the discrete criterion compares with its
continuous counterpart.

\begin{remark}[Comparison with the continuous problem]
\label{rem:cont} The harmonic duality of Theorem~\ref{thm:dual} has a
classical analogue. Let $U\Subset W$ be bounded smooth domains in $\mathbb{R}%
^{d}$, let $u\in C^{2}(\overline{U})$ vanish on $\partial U$, and let $%
f=-\Delta u$ in $U$. For the zero extensions $u^{0},f^{0}$ to $W$,
integration by parts gives, for $\phi \in C_{c}^{\infty }(W)$,
\begin{equation}
\langle -\Delta u^{0}-f^{0},\phi \rangle =\int_{\partial U}\partial _{\nu
}u\,\phi \,dS.  \label{eq:green}
\end{equation}%
Hence zero extension preserves the distributional equation exactly when the
normal derivative vanishes on the interface; Green's identity with a smooth
harmonic $h$ shows that this is equivalent to orthogonality of $f$ to
harmonic functions on $U$. The graph statement $\mathcal{L}_{SF}u=0$ and the
orthogonality in Theorem~\ref{thm:dual} are the discrete counterparts of
these two descriptions, and in both the elementary argument does not require
unique continuation. The graph version is finite-dimensional and computable,
with $\mathcal{L}_{SF}$ playing the role of the normal-derivative trace. To
the best of our knowledge, no result in this form---an interface criterion
with an exact error identity for zero extension of Dirichlet solutions---has
been reported on Riemannian manifolds, where the corresponding statements
are phrased through normal derivatives, harmonic orthogonality, or unique
continuation rather than through a finite interface operator. We do not
claim that the graph theorem implies the manifold theorem: a genuine
continuous limit would additionally require a graph construction, rescaling,
a topology of convergence, boundary control and estimates for the relevant
operators.
\end{remark}

\section{Conclusion}

\label{sec:concl}

The interface part of the graph Laplacian identifies the precise obstruction
to zero extension of a graph Dirichlet solution: the extension solves the
enlarged problem exactly when the interface block annihilates the solution,
and no condition is imposed on the exterior boundary of the enlarged domain.
The admissible sources form the orthogonal complement of the harmonic
extensions from the enlarged domain to the original set, and their space has
dimension equal to the number of vertices of the original set minus the rank
of the interface block; when the positive interface weights vary
independently, this dimension is generically equal to the number of vertices
minus the maximum matching size of the interface bipartite graph, while a
nonzero source of one sign is never admissible. When compatibility fails,
the Schur complement of the graph Laplacian on the enlarged domain with
respect to the original set gives the exact energy of the extension error
and two-sided spectral bounds, and the nearest admissible source is the
orthogonal projection onto the admissible space; exact cancellation can
disappear under an arbitrarily small change of the weights even while the
solution error tends to zero. The same interface reduction covers other
boundary conditions, including the Neumann problem, so the criterion is not
tied to the Dirichlet setting; it is the discrete counterpart of the
boundary-flux and harmonic-orthogonality descriptions of the continuous
problem, no result in this form has been reported on Riemannian manifolds,
and its higher-order extension requires boundary data compatible with the
lower-order interface traces.

\bigskip

\textbf{Data availability.} Data will be made available on request.

\textbf{Acknowledgments. }This work was supported by Doctoral guidance
program of Guangzhou Xinhua University. During the preparation of this work,
the authors used generative AI tools to improve the language clarity,
sentence structure, and overall readability of the manuscript. The
AI-assisted content was fully reviewed, revised, and optimized by the
authors. All scientific viewpoints, mechanism analysis, experimental
conclusions, and data interpretation are independently proposed and verified
by the authors. The authors take full responsibility for the accuracy and
integrity of all contents in this manuscript.


\begin{thebibliography}{99}
\bibitem{Bend} E. Bendito, A. Carmona, A. M. Encinas, Solving boundary value
problems on networks using equilibrium measures, J. Funct. Anal. 171 (2000)
155--176.

\bibitem{Cai} Y. Cai, S. Zhou, Zero extension for Poisson's equation, Sci.
China Math. 63 (2020) 721--732.

\bibitem{Chen} G.-H. Chen, K. Li, System and method for controlling
zero-count errors in computed tomography, U.S. Patent Application
2024/0177375 A1 (2024).

\bibitem{Dorfler} F. D\"orfler, F. Bullo, Kron reduction of graphs with
applications to electrical networks, IEEE Trans. Circuits Syst. I 60 (2013)
150--163.

\bibitem{Davis} T. A. Davis, S. Rajamanickam, W. M. Sid-Lakhdar, A survey of
direct methods for sparse linear systems, Acta Numer. 25 (2016) 383--566.

\bibitem{Huo} Z. Huo, F. Wiesinger, J. D. Kaggie, M. J. Graves, Low-rank
iterative infilling for zero echo-time (ZTE) imaging, Magn. Reson. Med. 93
(2025) 1149--1162.

\bibitem{Lian} Y. Lian, Z. Liu, J. Wang, S. Lu, Spatial-frequency aware
zero-centric residual unfolding network for MRI reconstruction, Magn. Reson.
Imaging 117 (2025) 110334.

\bibitem{mon} A. Mongia, A. Majumdar, Drug-target interaction prediction
using multi graph regularized nuclear norm minimization, PLoS ONE 15 (2020)
e0226484.

\bibitem{Shi} K. Shi, M. Burger, Continuum limit of $p$-biharmonic equations
on graphs, SIAM J. Math. Anal. 57 (2025) 3784--3805.

\bibitem{Wangm} M. N. Wang, Y. Li, L. L. Lei, D. W. Ding, X. J. Xie,
Combining non-negative matrix factorization with graph Laplacian
regularization for predicting drug-miRNA associations based on multi-source
information fusion, Front. Pharmacol. 14 (2023) 1132012.

\bibitem{Xuz} S. P. Xu, S. L. Zhou, Zero extension for the biharmonic
equation, Acta Math. Sin. (Engl. Ser.) 34 (2018) 1549--1562.

\bibitem{xuy} S. Xu, C. Yu, Zero extension for the Dirichlet problem of the
biharmonic equation, Mathematics 13 (2025) 1774.
\end{thebibliography}
\end{document}